\documentclass[12pt]{article}
\usepackage{epsfig}
\usepackage{ulem}
\usepackage{verbatim}
\usepackage{graphicx}
\usepackage{epstopdf}
\usepackage{color}
\usepackage{amssymb}
\usepackage{amsmath}
\usepackage{dsfont}
\usepackage{psfrag}

\usepackage{graphicx}
\usepackage{epstopdf}
\usepackage{color}
\usepackage{amssymb}
\usepackage{amsmath}
\usepackage{dsfont}
\usepackage{psfrag}

\graphicspath{{examples/}}

\newtheorem{lemma}{Lemma}
\newtheorem{corollary}{Corollary}
\newtheorem{proposition}{Proposition}

\newcommand{\prf}[1]{{\bf Proof} \, #1 \hfill $\blacksquare$}

\newcommand{\R}{\mathds{R}}
\newcommand{\C}{\mathds{C}}

\newcommand{\mat}[1]{\begin{bmatrix} #1 \end{bmatrix}}

\newcommand{\E}[1]{{\bf E}\left[#1\right]}

\newcommand{\ra}{\rightarrow}

\newcounter{l1}
\newcommand{\barablist}{\begin{list}{\arabic{l1}}{\usecounter{l1}}}

\title{\LARGE \bf
On Base Station Localization for State Estimation over Lossy
  Networks}

\author{Ufuk Topcu, Kenneth Hsu, and Kameshwar Poolla % <-this % stops a space
\thanks{The authors contributed equally to this work.}
\thanks{This work was supported in part by the NSF under Grant ECS
  03-02554.}% <-this % stops a space
\thanks{Authors are with the Department of Mechanical Engineering at the University of California at Berkeley.}
       }

\begin{document}

\maketitle
\thispagestyle{empty}
\pagestyle{empty}

%%%%%%%%%%%%%%%%%%%%%%%%%%%%%%%%%%%%%%%%%%%%%%%%%%%%%%%%%%%%%%%%%%%%%%
% Abstract
%%%%%%%%%%%%%%%%%%%%%%%%%%%%%%%%%%%%%%%%%%%%%%%%%%%%%%%%%%%%%%%%%%%%%%

\begin{abstract}
We consider a state estimation problem where observations
are made by multiple sensors.  These observations are communicated
over a lossy wireless network to a central base station that computes
estimates via a Kalman filter.  The goal is to determine the optimal
location of the base station under a certain class of packet loss
probability models.  It is shown in the two sensor case that the base
station is optimally located at one of the sensor locations.
Empirical evidence suggests that the result holds in some generality.

\end{abstract}

%%%%%%%%%%%%%%%%%%%%%%%%%%%%%%%%%%%%%%%%%%%%%%%%%%%%%%%%%%%%%%%%%%
% Body
%%%%%%%%%%%%%%%%%%%%%%%%%%%%%%%%%%%%%%%%%%%%%%%%%%%%%%%%%%%%%%%%%%

%%%%%%%%%%%%%%%%%%%%%%%%%%%%%%%%%%%%%%%%%%%%%%%%%%%%%%%%%%%%%%%%%%%%%%
% Introduction
%%%%%%%%%%%%%%%%%%%%%%%%%%%%%%%%%%%%%%%%%%%%%%%%%%%%%%%%%%%%%%%%%%%%%%

\section{Introduction}
\label{sec:Introduction}

The recent confluence of low-cost sensing, communication, and
computation technologies has led to much interest in the development
of wireless sensor networks for estimation and control.  These
wireless sensor networks provide for an economical means of
extracting greater performance and efficiency in a variety of
applications, such as manufacturing and chemical plant processes
\cite{Moyne}, indoor climate control (HVAC) \cite{Pakzad}, environment
monitoring \cite{Szewczyk}, electrical power distribution
\cite{Chong}, and automatic traffic flow \cite{Giridhar}.

There are many advantages that wireless sensor networks hold over
their wired counterparts.  For example, the fact that communication is
performed wirelessly eliminates any physical connection between the
various nodes in the network.  This enables wireless sensor networks
to be implemented without any major infrastructure overhauls.  This
also allows for the design of networks that can gracefully accept
dynamic changes to its structure, such as the addition or subtraction
of nodes, or in the case of mobile sensor networks, changing
communication topologies.

While the physical disconnection within the network allows for many
practical and performance related advantages, it also introduces
many fundamental issues (e.g. loss/delay of information,
limits on communication bandwidth, power constraints
\cite{Culler}) that are of lesser concern in wired configurations.
Issues such as packet loss/delay are crucial in estimation and control
problems where information flow is assumed to be uninterrupted, and
much work has been done to examine performance in the presence of
these complications \cite{Sinopoli,Schenato,Hespanha}.

In this paper, we consider the problem where (noisy) observations of
a dynamic process are made by multiple sensors at fixed locations.
These sensors communicate over a lossy wireless network to a base
station that computes the optimal state estimate by processing the
data via a Kalman filter.  Since the physical distance between the
sensors and the base station affects the packet loss probability
(which in turn affects the quality of the estimates), our goal is to
determine the optimal position of the base station.  A similar problem
of locating control logic has been studied in \cite{Kumar}.

There has been much recent work on decentralized estimation problems
where no single node is burdened with the majority of the
computational tasks \cite{Grime,Spanos}.  However, implementation of
these ideas require ``smart'' sensors that have local computational
abilities.  While these approaches have been shown to possess various
performance advantages, they are not always realizable. In order to
examine more fundamental issues regarding the implementation of
wireless sensor networks, we restrict our attention to the case where
stationary sensors merely \emph{sense} and \emph{communicate}, and
possess no computational power.  We leave the more general problem for
future work.

The remainder of the paper is organized as follows.  In Section
\ref{sec:Preliminaries}, we give a brief overview of Kalman
filtering.  Section \ref{sec:Base Station Placement} contains our
problem formulation and main results.  Some discussion regarding more
general cases can be found in Section \ref{sec:General Cases}.

% \input{notation.tex}
%%%%%%%%%%%%%%%%%%%%%%%%%%%%%%%%%%%%%%%%%%%%%%%%%%%%%%%%%%%%%%%%%%%%%%
% Preliminaries
%%%%%%%%%%%%%%%%%%%%%%%%%%%%%%%%%%%%%%%%%%%%%%%%%%%%%%%%%%%%%%%%%%%%%%

\section{Preliminaries}
\label{sec:Preliminaries}

%%%%%%%%%%%%%%%%%%%%%%%%%%%%%%%%%%%%%%%%%%%%%%%%%%%%%%%%%%%%%%%%%%%%%%
% Kalman Filtering
%%%%%%%%%%%%%%%%%%%%%%%%%%%%%%%%%%%%%%%%%%%%%%%%%%%%%%%%%%%%%%%%%%%%%%

\subsection{Kalman Filtering}
\label{sec:Kalman Filtering}

We begin by giving a brief overview of the canonical state
estimation problem.  Consider a process whose dynamics are modeled
by the (possibly time-varying) state space equation
\begin{equation}
\label{eq:process dynamics}
x_{k + 1} = A x_k + w_k,
\end{equation}
where $x_k \in \R^n$ is the state vector at time $k$, and $w$ is a
sequence of independent Gaussian random vectors with zero mean and
covariance matrix $Q$. The initial condition $x_0$ is unknown and
modeled as a Gaussian random variable with zero mean and covariance
matrix $P_0$.

We monitor this process via measurements of the form
\[
y_k = C x_k + v_k,
\]
where $C$ is possibly time-varying, and $v$ is a sequence of
independent Gaussian random vectors with zero mean and covariance
matrix $R$.  The sequences $w$ and $v$ are assumed to be
uncorrelated. The minimum variance estimate of the sequence $x$ is
then given recursively by the Kalman filter \cite{Kalman},
\begin{subequations}
\label{eq:Kalman Filter}
\begin{eqnarray}
K_k & = & A P_k C^* (C P_k C^* + R)^{-1} \\
\hat{x}_{k + 1} & = & A \hat{x}_k + K_k (y_k - C \hat{x}_k) \\
P_{k + 1} & = & A P_k A^* + Q - K_k C P_k A^*.
\end{eqnarray}
\end{subequations}
Here, $M^*$ denotes the complex conjugate transpose of a matrix $M$.

Due to the increasing availability of high quality, low-cost
sensors, we are often confronted with the situation where we have
access to measurements from multiple sensors.   That is, suppose
that we monitor the process \eqref{eq:process dynamics} via $N$
sensors, whose measurements are represented as
\begin{eqnarray*}
y_k^1 & = & C_1 x_k + v_k^1 \\
& \vdots & \\
y_k^N & = & C_N x_k + v_k^N.
\end{eqnarray*}
Here, $v^1, \dots, v^2$ are sequences of independent Gaussian random
vectors with zero mean and covariance matrices $R_1, \dots, R_N$,
respectively.  The sequences $v^1, \dots, v^N$ and $w$ are assumed to
be uncorrelated.  The resulting state estimates are again given by
\eqref{eq:Kalman Filter}, where we now define
\begin{eqnarray*}
C & = & \mbox{blkcol}(C_j), \ \ j = 1, \dots, N \\
R & = & \mbox{blkdiag}(R_j), \ \ j = 1, \dots, N.
\end{eqnarray*}
Here, $\mbox{blkcol}$ denotes the vertical stacking of its
arguments, and $\mbox{blkdiag}$ denotes a block diagonal matrix of
its arguments.

%%%%%%%%%%%%%%%%%%%%%%%%%%%%%%%%%%%%%%%%%%%%%%%%%%%%%%%%%%%%%%%%%%%%%%
% Kalman Filtering with Packet Losses
%%%%%%%%%%%%%%%%%%%%%%%%%%%%%%%%%%%%%%%%%%%%%%%%%%%%%%%%%%%%%%%%%%%%%%

\subsection{Kalman Filtering with Packet Losses}
\label{sec:Kalman Filter with Packet Losses}

We now discuss the minimum variance state estimation problem in the
presence of packet losses.  Let $\Omega$ be the power set of $\{1,
\dots, N\}$ and let each subset $\omega \in \Omega$ be given the
standard simple order.  Since the Kalman filtering equations
\eqref{eq:Kalman Filter} yield the minimum variance state estimate
for linear time-varying systems, we obtain the following optimal
filter for state estimation with packet losses
\[
\left. \begin{array}{l} K_k = A P_k F_{\omega_k}^* (F_{\omega_k} P_k
F_{\omega_k}^* + G_{\omega_k})^{-1} \\
\hat{x}_{k + 1} = A \hat{x}_k + K_k (y_k - F_{\omega_k}\hat{x}_k)  \\
P_{k + 1} = A P_k A^* + Q - K_k F_{\omega_k} P_k A^*
\end{array} \right\} \mbox{ for } \omega_k \in \Omega \backslash \phi,
\]
and
\[
\left. \begin{array}{l}
\hat{x}_{k + 1} = A \hat{x}_k \\
P_{k + 1} = A P_k A^* + Q
\end{array} \right\} \mbox{ for } \omega_k = \phi,
\]
 where
\begin{eqnarray*}
F_{\omega_k} & = & \mbox{blkcol}(C_j), \ \ j \in \omega_k \\
G_{\omega_k} & = & \mbox{blkdiag}(R_j), \ \ j \in \omega_k.
\end{eqnarray*}
For a similar treatment, see \cite{Liu,Hounkpevi}.

%%%%%%%%%%%%%%%%%%%%%%%%%%%%%%%%%%%%%%%%%%%%%%%%%%%%%%%%%%%%%%%%%%%%%%
% Base Station Placement
%%%%%%%%%%%%%%%%%%%%%%%%%%%%%%%%%%%%%%%%%%%%%%%%%%%%%%%%%%%%%%%%%%%%%%

\section{Base Station Placement}
\label{sec:Base Station Placement}

%%%%%%%%%%%%%%%%%%%%%%%%%%%%%%%%%%%%%%%%%%%%%%%%%%%%%%%%%%%%%%%%%%%%%%
% Problem Formulation
%%%%%%%%%%%%%%%%%%%%%%%%%%%%%%%%%%%%%%%%%%%%%%%%%%%%%%%%%%%%%%%%%%%%%%

\subsection{Problem Formulation}
\label{sec:Problem Formulation}

We will model the packet arrival process for each sensor as a sequence
of independent binomial random variables.  More specifically, for $i =
1, \dots, N$, let $\lambda_i$ be the probability that the observation
from sensor $i$ is received.  Since the sensor locations are fixed, a
given base station location gives rise to distances $d_1, \dots, d_N$.
We will assume that the probability of packet arrival decreases as the
distance $d_i$ between a sensor and the base station increases. 

Note that the sequences $(F_{\omega})_k$ and $(G_{\omega})_k$ both
depend on whether packets have been received or lost, and hence the
sequence $(P_k)_k$ is itself random.  We will then aim to minimize
(with respect to the Loewner ordering \cite{Horn}, i.e., $X \succeq Y$ if
$X - Y$ is positive semi-definite). the expected estimation error
covariance.  That is, we consider the problem 
\[
\min \E{P_{k + 1} | P_k},
\]
where the minimization is performed over the set of possible base
station locations.  If we define for each nonempty $\neq \omega \in
\Omega$
\begin{eqnarray*}
\alpha_{\omega} & = & \prod_{j \in \omega} \lambda_j \prod_{j \in \{1,
\dots, N\} \backslash \omega} (1 - \lambda_j) \\
K_k & = & A P_k F_{\omega}^* (F_{\omega} P_k F_{\omega}^* +
G_{\omega})^{-1},
\end{eqnarray*}
and $\alpha_\omega = 0$ for $\omega = \phi$, we can write
\[
\E{P_{k + 1} | P_k} = A P_k A^* + Q - \sum_{\omega \in \Omega}
\alpha_\omega K_k F_\omega P_k A^*.
\]
As stated, the above optimization problem is a difficult nonlinear
programming problem that can be attacked by the appropriate solvers.
Since it is not feasible for the location of the base station to
change over time, we will insist that the solution to the above
problem minimizes $\E{P_{k + 1} | P_k}$ for any $P_k$.  Consequently,
we shall hereafter restrict our attention to an illuminating instance
of the above problem.

%%%%%%%%%%%%%%%%%%%%%%%%%%%%%%%%%%%%%%%%%%%%%%%%%%%%%%%%%%%%%%%%%%%%%%
% The 2 Sensor Case
%%%%%%%%%%%%%%%%%%%%%%%%%%%%%%%%%%%%%%%%%%%%%%%%%%%%%%%%%%%%%%%%%%%%%%

\subsection{The 2 Sensor Case}
\label{sec:The 2 Sensor Case}

In this section, we consider the case where measurements are obtained
from $2$ homogeneous sensors, $s_1$ and $s_2$.  Without loss of
generality, these sensors are located at $\xi_1 = 0$ and $\xi_2 = 1$
on the real line.  The goal is to find the best location $d \in [0,1]$
for the base station.  Note that optimal solution must lie in the
interval $[0,1]$.  The situation $d = 0$ (or $d = 1)$ corresponds to
the base station being physically wired to one of the sensors, and
wireless communication is performed only by the other sensor.

Although there may be many factors that influence the packet loss
probability, we shall consider only the effects due to the physical
distance between the communicating nodes.  Since the sensors have
identical broadcasting capabilities, it is then natural to model the
probability of packet loss in the following manner.   Let $f : [0,1]
\ra [0,1]$ be convex and decreasing, with $f(0) = 1$.  When the base
station is located at some $d \in [0,1]$, the probability that a
packet is received from sensor $s_1$ is $f(d)$, and the probability
that a packet is received from sensor $s_2$ is $f(1 - d)$.  The
constraint that $f(0) = 1$ captures the implicit assumption that no
packet is lost when a sensor is affixed to the base station and
communication is not performed wirelessly.

%%%%%%%%%%%%%%%%%%%%%%%%%%%%%%%%%%%%%%%%%%%%%%%%%%%%%%%%%%%%%%%%%%%%%%
% Vector-valued Measurements, Same Covariances
%%%%%%%%%%%%%%%%%%%%%%%%%%%%%%%%%%%%%%%%%%%%%%%%%%%%%%%%%%%%%%%%%%%%%%

\subsubsection{Vector-valued Measurements, Same Covariances}
\label{sec:Vector-valued Measurements, Same Covariances}

We will first restrict our attention to the case where 
\begin{eqnarray*}
C_1 & = & C_2 = C \in \R^{m \times n} \\
R_1 & = & R_2 = R \in \R^{m \times m}.
\end{eqnarray*}
If we model the probability of packet loss as described above, we can
examine the estimation error covariance as a function of $d$.  For
ease of notation, let us define
\begin{eqnarray*}
M_1 & = & A P_k C^* (C P_k C^* + R)^{-1} C P_k A^* \\
M_2 & = & A P_k C^* (C P_k C^* + R/2)^{-1} C P_k A^* \\
S & = & A P_k A^* + Q.
\end{eqnarray*}
The estimation error covariance $P_{k + 1}$ then assumes values as
follows:
\[
\begin{array}{|c|c|} \hline
P_{k + 1} & \mbox{Probability} \\ \hline
S - M_2 &  \ \ f(d) f(1 - d) \\ \hline
\ \ S - M_1 \ \ &  \ \ f(d) - 2 f(d) f(1 - d)) + f(1 - d) \ \ \\ \hline
S &  (1 - f(d)) (1 - f(1 - d)) \\ \hline
\end{array}
\]
% \[
% \begin{array}{ccc}
% S - M_2 & \mbox{ w.p. } & f(d) f(1 - d) \\
% S - M_1 & \mbox{ w.p. } & f(d) - 2 f(d) f(1 -
% d)) + f(1 - d) \\
% S & \mbox{ w.p. } & (1 - f(d)) (1 - f(1 - d))
% \end{array}
% \]
Note that in the case where packets from both sensors are lost, the
resulting covariance corresponds to that of propagating the previous
state estimate.  The expected error covariance can then be computed as
\begin{eqnarray*}
\label{eq:cost function}
\E{P_{k + 1} | P_k}
& = & S - f(d) f(1 - d) M_1 \\
& & \hspace{0.05 in} - (f(d) - 2 f(d) f(1 - d) + f(1 - d)) M_2.
\end{eqnarray*}

We have the following result.
\begin{proposition}
\label{thm:2 sensors}
Let $C_1 = C_2 = C \in \R^{m \times n}$ and $R_1 = R_2 = R \in \R^{m
\times m}$.  Let $f:[0,1] \ra [0,1]$ be twice differentiable, convex,
and decreasing, with $f(0) = 1$.  Then, $\E{P_{k + 1} | P_k}$ is
matrix concave in $d$ on $[0,1]$.
\end{proposition}

\prf{
Define
\begin{eqnarray*}
J_{k + 1}(d) & = & \E{P_{k + 1} | P_k} \\
& = & A P_k A^* + Q - f(d) f(1 - d) M_2 \\
& & \hspace{0.15 in} - (f(d) - 2 f(d) f(1 - d) + f(1 - d)) M_1.
\end{eqnarray*}
% Taking the derivative with respect to $d$, we have
% \begin{eqnarray*}
% J_{k + 1}'(d) & = & (f'(d) f(1 - d) - f(d) f'(1 - d)) (2 M_1 - M_2) \\
% & & \hspace{0.15 in} - (f'(d) + f'(1 - d)) M_1.
% \end{eqnarray*}
Taking the second derivative yields
\begin{eqnarray*}
J_{k + 1}''(d) & = &
% (f''(d) f(1 - d) - 2f'(d)f'(1 - d) \\
% & & \hspace{0.02 in} + f(d) f''(1 - d)) (2M_1 - M_2) \\
% & & \hspace{0.02 in} - (f''(d) + f''(1 - d))M_1 \\
f''(d) (f(1 - d) - 1) M_1 \\
& & \hspace{-0.03 in} + f''(1 - d) (f(d) - 1) M_1 \\
& & \hspace{-0.03 in} - 2 f'(d) f'(1 - d) (2 M_1 - M_2) \\
& & \hspace{-0.03 in} + (f''(d) f(1 - d) + f(d) f''(1 - d)) (M_1 - M_2).
\end{eqnarray*}
Since $f$ is convex and decreasing with $f(0) = 1$, we have $f'(x)
\leq 0$ and $f''(x) \geq 0$ for $x \in [0,1]$.  Consequently, $J''(d)
\preceq 0$, and hence $\E{P_{k + 1} | P_k}$ is matrix concave on
$[0,1]$ (see \cite{Boyd}, p.110). 
}

\begin{corollary}
\label{cor:2 sensors}
Let $C_1 = C_2 = C \in \R^{m \times n}$ and $R_1 = R_2 = R \in \R^{m
\times m}$.  Let $f:[0,1] \ra [0,1]$ be twice differentiable, convex,
and decreasing, with $f(0) = 1$.  Then, the base station is optimally
located at one of the sensor positions.
\end{corollary}

%%%%%%%%%%%%%%%%%%%%%%%%%%%%%%%%%%%%%%%%%%%%%%%%%%%%%%%%%%%%%%%%%%%%%%
% Scalar-valued Measurements, Same Covariances
%%%%%%%%%%%%%%%%%%%%%%%%%%%%%%%%%%%%%%%%%%%%%%%%%%%%%%%%%%%%%%%%%%%%%%

\subsubsection{Scalar-valued Measurements, Same Covariances}
\label{sec:Scalar-valued Measurements, Same Covariances}

A similar result may be obtained in the case where $C_1 = C_2  = C \in
\R^{1 \times n}$ and the noise covariances $R_1$ and $R_2$ are not
necessarily equal.  We will need the following preliminary result.
\begin{lemma}
\label{lem:M_0 - M_1 - M_2}
Let $C_1 = C_2 = C \in \R^{1 \times n}$ and define
\begin{eqnarray*}
T & = & C P_k C^* \\
M_0 & = & AP\mat{C\\C}^* \mat{T + R_1 & T \\ T & T + R_2
}^{-1} \mat{C\\C} PA^* \\
M_1 & = & APC^* (T + R_1)^{-1}C PA^* \\
M_2 & = & APC^* (T + R_2)^{-1}C PA^*.
\end{eqnarray*}
Then, $M_0 - M_1 - M_2 \succeq 0$.
\end{lemma}

\prf{
%Let $\alpha$, $\beta$ and $\gamma$ be defined as
%\[
%\mat{\alpha& \beta \\ \beta & \gamma} = \mat{T + R_1 & T \\
%T & T + R_2 }
%\]
%and let
%\[
%U = \mat{C \\ C} P_k A^*.
%\]
Let us define
\[
\mat{\alpha& \beta \\ \beta & \gamma} = \mat{T + R_1 & T \\
T & T + R_2 }, \ \mbox{ and } \
U = \mat{C \\ C} P_k A^*.
\]
Then, $M_0 - M_1 - M_2$ can be written as
\begin{eqnarray*}
M_0 - & M_1 & - M_2 \\
& = & U^* \left( \mat{\alpha & \beta \\ \beta & \gamma }^{-1} -
  \mat{\alpha & 0 \\ 0 & \gamma}^{-1} \right) U \\
& = & U^* \mat{\frac{1}{\alpha - \beta^2 \gamma^{-1}} -\alpha^{-1}
& - \frac{\beta \gamma^{-1}}{\alpha - \beta^2 \gamma^{-1}} \\
  -\frac{\gamma^{-1} \beta}{\alpha - \beta^2 \gamma^{-1}} &
  \frac{\gamma^{-2}\beta^2}{\alpha - \beta^2 \gamma^{-1}}} U \\
& = & U^* \mat{\frac{1}{\alpha - \beta^2 \gamma^{-1}} \\
  -\frac{\gamma^{-1} \beta}{\alpha - \beta^2 \gamma^{-1}}}
  \mat{\frac{1}{\alpha - \beta^2 \gamma^{-1}}  &
  -\frac{\beta\gamma^{-1}}{\alpha - \beta^2 \gamma^{-1}}} U \\
& & \hspace{0.15 in} - U^* \mat{-\alpha^{-1} & 0 \\ 0 & 0} U \\
& = & A P_k C^* \left( \frac{(1 -\gamma^{-1} \beta)^2}{\alpha -
\beta^2 \gamma^{-1}} - \alpha^{-1} \right) C P_k A^*,
\end{eqnarray*}
where the second equality follows from the matrix inversion lemma
\cite{Horn}. Now, let
\begin{eqnarray*}
H & = & \frac{(1 -\gamma^{-1} \beta)^2}{\alpha - \beta^2 \gamma^{-1}}
- \alpha^{-1} \\
& = & \frac{(1-\gamma^{-1}\beta)^2 - \alpha^{-1} (\alpha -\beta^2
\gamma^{-1})}{\alpha -\beta^2 \gamma^{-1}}.
\end{eqnarray*}
The desired result follows from noticing that the denominator is
positive and the numerator satisfies
\begin{eqnarray*}
& & (1 - \gamma^{-1} \beta)^2 - \alpha^{-1} (\alpha - \beta^2
\gamma^{-1}) \\
& & \hspace{0.35 in} = -\beta \gamma^{-1} \left(2 - \frac{C P_k C^*}{C
P_k C^* + R_2} - \frac{C P_k C^*}{C P_k C^* + R_1} \right) \\
& & \hspace{0.35 in} \leq 0.
\end{eqnarray*}
}

\begin{proposition}
Let $C_1 = C_2 = C \in \R^{1 \times n}$ and suppose that $R_1$ and
$R_2$ are not necessarily equal.  Let $f:[0,1] \ra [0,1]$ be twice
differentiable, convex, and decreasing, with $f(0) = 1$.  Then,
$\E{P_{k + 1} | P_k}$ is matrix concave in $d$ on $[0,1]$.
\end{proposition}

\prf{
Define
\begin{eqnarray*}
T & = & C P_k C^* \\
M_0 & = & AP\mat{C\\C}^* \mat{T + R_1 & T \\ T & T + R_2
}^{-1} \mat{C\\C} PA^* \\
M_1 & = & APC^* (T + R_1)^{-1}C PA^* \\
M_2 & = & APC^* (T + R_2)^{-1}C PA^*.
\end{eqnarray*}
We can then write
\begin{eqnarray*}
J_{k + 1}(d) & = & \E{P_{k + 1} | P_k} \\
& = & A P_k A^* + Q \\
& & \hspace{0.15 in} - f(d) f(1 - d) M_0 - f(d) (1 - f(1 - d)) M_1 \\
& & \hspace{0.15 in} - f(1 - d) (1 - f(d)) M_2 \\
\end{eqnarray*}
% Taking the derivative with respect to $d$, we have
% \begin{eqnarray*}
% J_{k + 1}'(d) & = & (f'(d) f(1 - d) \\
% & & \hspace{0.15 in} - f(d) f'(1 - d)) (M_1 + M_2 - M_0) \\
% & & \hspace{0.15 in} - f'(d) M_1 - f'(1 - d)) M_2.
% \end{eqnarray*}
Taking the second derivative yields
\begin{eqnarray*}
J_{k + 1}''(d) & = &
% (f''(d) f(1 - d) - 2f'(d)f'(1 - d) \\
% & & \hspace{0.02 in} + f(d) f''(1 - d)) (2M_1 - M_2) \\
% & & \hspace{0.02 in} - (f''(d) + f''(1 - d))M_1 \\
f''(d) (f(1 - d) - 1) M_1 \\
& & \hspace{-0.03 in} + f''(1 - d) (f(d) - 1) M_2 \\
& & \hspace{-0.03 in} - 2 f'(d) f'(1 - d) (M_1 + M_2 - M_0) \\
& & \hspace{-0.03 in} + f''(d) f(1 - d) (M_2 - M_0) \\
& & \hspace{-0.03 in} + f(d) f''(1 - d) (M_1 - M_0).
\end{eqnarray*}
Since $f$ is convex and decreasing with $f(0) = 1$, we have $f'(x)
\leq 0$ and $f''(x) \geq 0$ for $x \in [0,1]$.  From Lemma
\ref{lem:M_0 - M_1 - M_2} we have $M_1 + M_2 - M_0 \preceq 0$.  As a
result, $ J''(d) \preceq 0$, so $\E{P_{k + 1} | P_k}$ is matrix
concave on $[0,1]$ (see \cite{Boyd}, p.110). 
}

\begin{corollary}
Let $C_1 = C_2 = C \in \R^{1 \times n}$ and suppose that $R_1$ and
$R_2$ are not necessarily equal.  Let $f : [0,1] \ra [0,1]$ be twice
differentiable, convex, and decreasing, with $f(0) = 1$.  Then, the
base station is optimally located at one of the sensor positions.  If
$R_1 \succeq R_2$, then the base station is optimally located at the
second sensor.
\end{corollary}

%%%%%%%%%%%%%%%%%%%%%%%%%%%%%%%%%%%%%%%%%%%%%%%%%%%%%%%%%%%%%%%%%%%%%%
% An Example
%%%%%%%%%%%%%%%%%%%%%%%%%%%%%%%%%%%%%%%%%%%%%%%%%%%%%%%%%%%%%%%%%%%%%%

\subsubsection{An Example}
\label{sec:An Example}

We now offer an illustrative example to complement our main
results.

\begin{figure}[h]
\centering 
% This file is generated by the MATLAB m-file laprint.m. It can be included
% into LaTeX documents using the packages graphicx, color and psfrag.
% It is accompanied by a postscript file. A sample LaTeX file is:
%    \documentclass{article}\usepackage{graphicx,color,psfrag}
%    \begin{document}\input{ex1}\end{document}
% See http://www.mathworks.de/matlabcentral/fileexchange/loadFile.do?objectId=4638
% for recent versions of laprint.m.
%
% created by:           LaPrint version 3.16 (13.9.2004)
% created on:           31-May-2007 13:40:39
% eps bounding box:     10 cm x 7.5 cm
% comment:              
%
\begin{psfrags}%
\psfragscanon%
%
% text strings:
\psfrag{s01}[t][t]{\color[rgb]{0,0,0}\setlength{\tabcolsep}{0pt}\begin{tabular}{c}Time Index\end{tabular}}%
\psfrag{s02}[b][b]{\color[rgb]{0,0,0}\setlength{\tabcolsep}{0pt}\begin{tabular}{c}${\rm Trace}(P_k)$\end{tabular}}%
%
% xticklabels:
\psfrag{x01}[t][t]{0}%
\psfrag{x02}[t][t]{10}%
\psfrag{x03}[t][t]{20}%
\psfrag{x04}[t][t]{30}%
\psfrag{x05}[t][t]{40}%
\psfrag{x06}[t][t]{50}%
\psfrag{x07}[t][t]{60}%
\psfrag{x08}[t][t]{70}%
\psfrag{x09}[t][t]{80}%
\psfrag{x10}[t][t]{90}%
\psfrag{x11}[t][t]{100}%
%
% yticklabels:
\psfrag{v01}[r][r]{2}%
\psfrag{v02}[r][r]{2.2}%
\psfrag{v03}[r][r]{2.4}%
\psfrag{v04}[r][r]{2.6}%
\psfrag{v05}[r][r]{2.8}%
\psfrag{v06}[r][r]{3}%
\psfrag{v07}[r][r]{3.2}%
\psfrag{v08}[r][r]{3.4}%
\psfrag{v09}[r][r]{3.6}%
\psfrag{v10}[r][r]{3.8}%
\psfrag{v11}[r][r]{4}%
%
% Figure:
\resizebox{8cm}{!}{\includegraphics{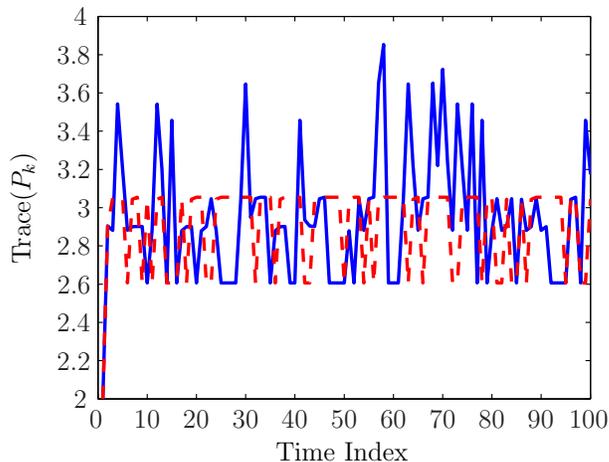}}%
\end{psfrags}%
%
% End ex1.tex

% \includegraphics[width = 2.5 in]{examples/ex1.eps}
\caption{Covariances for $d = 0$ (dashed) and $d = 0.5$ (solid).} 
\label{fig:ex1}
\end{figure}

Suppose that there are two sensors, located respectively at the points
$0$ and $1$ on the real line.  Let us model our packet arrival
probability function as
\[
f(d) = e^{-d}.
\]
Note that $f$ is convex and decreasing on $[0,1]$ with $f(0) = 1$.
The process being monitored is described by the state space equations
in \eqref{eq:process dynamics}, where all the eigenvalues of $A \in
\R^{2 \times 2}$ lie in the unit circle and $C \in \R^{1 \times 2}$.
The trace of the estimation error covariance matrices is plotted in
Figure \ref{fig:ex1} for the case where the base station is situated
at the point $0$ and the case where it is situated at the point $0.5$.

%%%%%%%%%%%%%%%%%%%%%%%%%%%%%%%%%%%%%%%%%%%%%%%%%%%%%%%%%%%%%%%%%%%%%%
% General Cases
%%%%%%%%%%%%%%%%%%%%%%%%%%%%%%%%%%%%%%%%%%%%%%%%%%%%%%%%%%%%%%%%%%%%%%

\section{General Cases}
\label{sec:General Cases}

We now offer some commentary on extending our main results.

\begin{figure}[h]
\centering
% This file is generated by the MATLAB m-file laprint.m. It can be included
% into LaTeX documents using the packages graphicx, color and psfrag.
% It is accompanied by a postscript file. A sample LaTeX file is:
%    \documentclass{article}\usepackage{graphicx,color,psfrag}
%    \begin{document}\input{ex2}\end{document}
% See http://www.mathworks.de/matlabcentral/fileexchange/loadFile.do?objectId=4638
% for recent versions of laprint.m.
%
% created by:           LaPrint version 3.16 (13.9.2004)
% created on:           31-May-2007 13:02:57
% eps bounding box:     10 cm x 7.5 cm
% comment:              
%
\begin{psfrags}%
\psfragscanon%
%
% text strings:
\psfrag{s01}[t][t]{\color[rgb]{0,0,0}\setlength{\tabcolsep}{0pt}\begin{tabular}{c}Location\end{tabular}}%
\psfrag{s02}[b][b]{\color[rgb]{0,0,0}\setlength{\tabcolsep}{0pt}\begin{tabular}{c}Expected error covariance\end{tabular}}%
%
% xticklabels:
\psfrag{x01}[t][t]{0}%
\psfrag{x02}[t][t]{0.1}%
\psfrag{x03}[t][t]{0.2}%
\psfrag{x04}[t][t]{0.3}%
\psfrag{x05}[t][t]{0.4}%
\psfrag{x06}[t][t]{0.5}%
\psfrag{x07}[t][t]{0.6}%
\psfrag{x08}[t][t]{0.7}%
\psfrag{x09}[t][t]{0.8}%
\psfrag{x10}[t][t]{0.9}%
\psfrag{x11}[t][t]{1}%
%
% yticklabels:
\psfrag{v01}[r][r]{1.22}%
\psfrag{v02}[r][r]{1.24}%
\psfrag{v03}[r][r]{1.26}%
\psfrag{v04}[r][r]{1.28}%
\psfrag{v05}[r][r]{1.3}%
\psfrag{v06}[r][r]{1.32}%
\psfrag{v07}[r][r]{1.34}%
%
% Figure:
\resizebox{8cm}{!}{\includegraphics{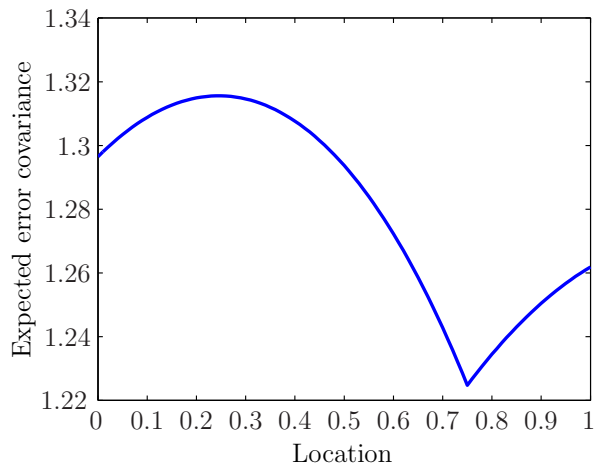}}%
\end{psfrags}%
%
% End ex2.tex

% \includegraphics[width = 2.5 in]{examples/ex2.eps}
\caption{Piecewise concavity for multiple colinear sensors.}
\label{fig:ex2}
\end{figure}

Although the previous results were confined to the case $C_1 = C_2$,
simulations suggest that they also hold in more generality.
Consider the case where $N$ sensors located colinearly are
available to make observations.  In general, we may have different
values for $C_1, \dots, C_N$ and different values for $R_1, \dots,
R_N$.  In this case, empirical evidence supports the extension of our
main result in the sense that the expected error covariance is
piecewise concave over the intervals defined by the location of the
sensors.  Figure \ref{fig:ex2} illustrates this situation for the case
$N = 3$.

\begin{figure}[h]
\centering
\input{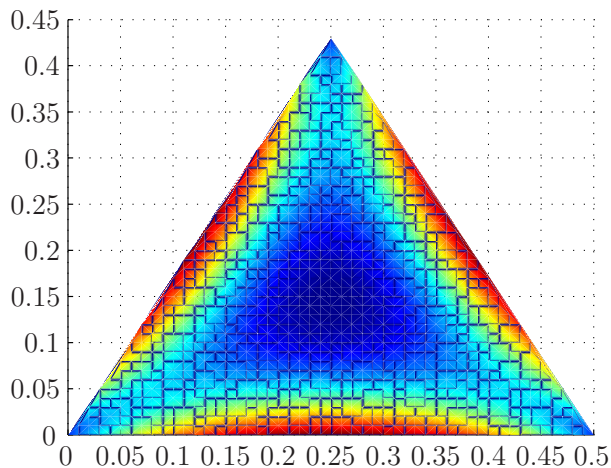}
\caption{Local minimum in the interior of the convex hull.}
\label{fig:ex3}
\end{figure}

In the case where the convex hull of the sensor locations has a
nonempty interior in $\R^2$, it can be shown that the optimal base
station location may lie at a point other than the sensor locations.
Consider 3 sensors, located at the vertices of an equilateral triangle
with side length $0.5$.  For $P_k = I$, $f(d) = e^{-d}$, $C_1 = C_2 =
C_3 = 0.6$, and unit noise covariances, the expected error covariance
(as a function of base station location) has a local minimum at the
centroid, as illustrated in Figure \ref{fig:ex3}.

%%%%%%%%%%%%%%%%%%%%%%%%%%%%%%%%%%%%%%%%%%%%%%%%%%%%%%%%%%%%%%%%%%%%%%
% Conclusion
%%%%%%%%%%%%%%%%%%%%%%%%%%%%%%%%%%%%%%%%%%%%%%%%%%%%%%%%%%%%%%%%%%%%%%

\section{Conclusion}
\label{sec:Conclusion}

We considered the problem of positioning a central computational node
within a wireless sensor network for the purpose of state estimation.
In the presence of potential packet losses, the optimal estimate is
given by a time-varying Kalman filter.  We then strived to determine
the optimal location of the base station so as to minimize the
expected conditional estimation error covariance.  It was shown for
the case of a 2-sensor configuration that under a certain class of
packet loss probability models, the base station is optimally located
at one of the sensor positions.

\bibliographystyle{plain}
\bibliography{references}

\begin{thebibliography}{10}

\bibitem{Boyd}
S.~Boyd and L.~Vandenberghe.
\newblock {\em Convex Optimization}.
\newblock Cambridge University Press, Cambridge, 2004.

\bibitem{Culler}
D.E. Culler, D.Estrin, and M.B. Srivastava.
\newblock Guest editors' introduction, overview of sensor networks.
\newblock {\em IEEE Computer}, 37(8):41--49, 2004.

\bibitem{Chong}
{C.Y. Chong and S.P. Kumar and B.A. Hamilton}.
\newblock Sensor networks, evolution, opportunities, and challenges.
\newblock {\em Proceedings of the IEEE}, 91:1247-- 1256, 2003.

\bibitem{Giridhar}
A.~Giridhar and P.R. Kumar.
\newblock Scheduling automated traffic on a network of roads.
\newblock {\em IEEE Transactions on Vehicular Technology}, 55:1467--1474, 2006.

\bibitem{Grime}
S.~Grime and H.R. Durrant-Whyte.
\newblock Data fusion in decentralized sensor networks.
\newblock {\em Control Engineering Practice}, 2:849--863, 1994.

\bibitem{Hespanha}
J.P. Hespanha, P.~Naghshtabrizi, and Y.~Xu.
\newblock A survey of recent results in networked control systems.
\newblock {\em Proceedings of the IEEE}, 95:138--162, 2007.

\bibitem{Horn}
R.A. Horn and C.R. Johnson.
\newblock {\em Matrix Analysis}.
\newblock Cambridge University Press, Cambridge, 1990.

\bibitem{Hounkpevi}
F.O. Hounkpevi and E.E. Yaz.
\newblock Robust minimum variance linear state estimators for multiple sensors
  with diffferent failure rates.
\newblock {\em Automatica}, 43:1274--1280, 2007.

\bibitem{Kalman}
R.E. Kalman.
\newblock A new approach to linear filtering and prediction problems.
\newblock {\em ASME Journal of Basic Engineering}, 82:35--45, 1960.

\bibitem{Liu}
X.~Liu and A.~Goldsmith.
\newblock Kalman filtering with partial observation losses.
\newblock {\em Proceedings of the IEEE Conference on Decision and Control},
  pages 4180--4186, 2004.

\bibitem{Moyne}
J.R. Moyne and D.M. Tilbury.
\newblock The emergence of industrial control networks for manufacturing
  control, diagnostics, and safety data.
\newblock {\em Proceedings of the IEEE}, 95:29--47, 2007.

\bibitem{Pakzad}
S.N. Pakzad, S.~Kim, G.L. Fenves, S.D. Glaser, D.E. Culler, and J.W. Demmel.
\newblock Multi-purpose wireless accelerometers for civil infrastructure
  monitoring.
\newblock {\em Proceedings of the 5th International Workshop on Structural
  Health Monitoring}, 2005.

\bibitem{Kumar}
C.L. Robinson and P.R. Kumar.
\newblock Control over networks of lossy links.
\newblock {\em IEEE Conference on Computer Communications}, 2007.

\bibitem{Schenato}
L.~Schenato, B.~Sinopoli, M.~Franceschetti, K.~Poolla, and S.~Sastry.
\newblock Foundations of control and estimation over lossy networks.
\newblock {\em IEEE Special Issue on Networked Control Systems}, 95:163--187,
  2007.

\bibitem{Sinopoli}
B.~Sinopoli, L.~Schenato, M.~Franceshetti, K.~Poolla, M.~Jordan, and S.~Sastry.
\newblock Kalman filtering with intermittent observations.
\newblock {\em IEEE Transactions on Automatic Control}, 49:1453--1469, 2004.

\bibitem{Spanos}
D.P. Spanos, R.~Olfati-Saber, and R.M. Murray.
\newblock Approximate distributed kalman filtering in sensor networks with
  quantifiable performance.
\newblock {\em Conference on Information Processing in Sensor Networks}, 2005.

\bibitem{Szewczyk}
R.~Szewczyk, J.~Polastre, A.~Mainwaring, J.~Anderson, and D.~Culler.
\newblock An analysis of a large scale habitat monitoring application.
\newblock {\em ACM Conference on Embedded Networked Sensor Systems}, pages
  214--226, 2004.

\end{thebibliography}

\end{document}